\documentclass[12pt]{amsart}
\usepackage{graphicx}
\usepackage[centertags]{amsmath}
\usepackage{amsfonts}
\usepackage{amssymb}
\usepackage{amsthm}
\usepackage{newlfont}
\newtheorem{thm}{Theorem}[section]

\newtheorem{prop}[thm]{Proposition}
\theoremstyle{definition}

\theoremstyle{remark}

\newcommand{\grad}{\textrm{grad}}

\title{On a class of special Riemannian manifold}
\bigskip
\author{Dimitar Razpopov}
\date{}
\begin{document}

\maketitle

\begin{abstract}
We consider a $4$-dimensional Riemannian manifold $M$ with a
metric $g$ and an affinor structure $q$. We note the local
coordinates of $g$ and $q$ are circulant matrices. Their first
orders are $(A, B, C, B)$, $A, B, C \in FM$ and $(0, 1, 0, 0)$,
respectively.

Let $\nabla$ be the connection of $g$. Then we
obtain:

1) \quad  $q^{4}=E;$\quad $g(qx, qy)=g(x,y)$,\quad $x,\ y\in \chi
M$,

2) \quad $\nabla q =0$ if and only if $
    \grad A=(\grad C)q^{2};\quad 2\grad B= (\grad C)(q+q^{3})$,

\end{abstract}

\Small{\textbf{Mathematics Subject Classification (2010)}: 53C15,
53B20}

\Small{\textbf{Keywords}: Riemannian metric, affinor structure,
sectional curvatures}
\thispagestyle{empty}
\section{Introduction}
The main purpose of the present paper is to find a class of
Riemannian manifolds which admits a circulant metric $g$, as well
as an additional circulant structure $q$, such that $q^{4}=id$,
and $q$ is a parallel structure with respect to the Riemannian
connection $\nabla$ of $g$.
\section{Preliminaries}
 We consider a $4$-dimensional Riemannian manifold $M$ with a
metric $g$ and an affinor structure $q$. We note the local
coordinates of $g$ and $q$ are circulant matrices. The analogous
case of a three dimensional Riemannian manifold has been discussed
in \cite{1}, \cite{2}, \cite{3}.

So let the metric $g$ have the coordinates:
\begin{equation}\label{f1}
    g_{ij}=\begin{pmatrix}
      A & B & C & B \\
      B & A & B & C\\
      C & B & A & B\\
      B & C & B & A\\
    \end{pmatrix},\qquad det g_{ij}=(A-C)^{2}((A+C)^{2}-4B^{2})
\end{equation}
in the local coordinate system $(x_{1}, x_{2}, x_{3}, x_{4})$, and
$A=A(p), B=B(p), C=C(p)$, where $p(x_{1}, x_{2}, x_{3}, x_{4})\in
F\subset R^{4}$.  Naturally, $A, B, C$ are smooth functions of a
point $p$. We suppose $A > C > B > 0$. These conditions imply the
heat minors of the matrix $g$ are positive, so the metric $g$ is
positively defined \cite{5}. The inverse matrix is
\begin{equation}\label{f2}
    g^{ij}=\frac{1}{D}\begin{pmatrix}
      \overline{A} & \overline{B} & \overline{C} & \overline{B} \\
      \overline{B} & \overline{A} & \overline{B} & \overline{C}\\
      \overline{C} & \overline{B} & \overline{A} &  \overline{B} \\
      \overline{B} & \overline{C} & \overline{B} & \overline{A}\\
    \end{pmatrix},\qquad D=(A-C)((A+C)^{2}-4B^{2}),
\end{equation}
where $\overline{A}=A(A+C)-2B^{2}$, $\overline{B}=B(C-A)$,
$\overline{C}=2B^{2}-C(A+C)$.

Further, let the local coordinates of $q$ be
\begin{equation}\label{f3}
    q_{i}^{.j}=\begin{pmatrix}
      0 & 1 & 0 & 0\\
      0 & 0 & 1 & 0\\
      0 & 0 & 0 & 1\\
      1 & 0 & 0 & 0\\
    \end{pmatrix}.
\end{equation}

We will use the notation $\Phi_{i}=\dfrac{\partial \Phi}{\partial
x^{i}}$ for every smooth function $\Phi$ defined in $F$.


\section{The condition for a Parallel structure}
\begin{thm} Let $M$ be a $4$-dimensional Riemannian manifold with
a metric $g$ and an affinor structure $q$ with local coordinates
(\ref{f1}), (\ref{f3}), respectively. Then we have

\begin{equation}\label{2.1}
    q^{4}=E;\quad q^{2}\neq \pm E
\end{equation}
\begin{equation}\label{2.2}
    g(qu, qv)=g(u, v),\quad u, v \in \chi M,
\end{equation}
where $E$ is the unit matrix.
\end{thm}

\begin{proof}
The conditions (\ref{2.1}) follows directly from (\ref{f3}).
\end{proof}

Now let $u=(u^{1},u^{2}, u^{3}, u^{4})$ and $v=(v^{1},v^{2},
v^{3}, v^{4})$ be two vectors in $\chi M$. Using (\ref{f1}) and
(\ref{f3}) we calculate that $$g(u, u)=
A((u^{1})^{2}+(u^{2})^{2}+(u^{3})^{2}+(u^{4})^{2})+2B(u^{1}u^{2}+u^{1}u^{4}+u^{2}u^{3}+u^{3}u^{4})+2C(u^{1}u^{3}+u^{2}u^{4}),$$ $g(qu, qu)=g(u, u)$.
It is easily to see $$g(q^{3}u, q^{3}u)=g(q^{2}u, q^{2}u)=g(qu, qu)=g(u, u),$$ $$g(q^{3}u, q^{3}v)=g(q^{2}u, q^{2}v)=g(qu, qv)=g(u, v).$$ 

\begin{prop}Let $u=(u^{1},u^{2}, u^{3}, u^{4})$ and $v=(v^{1},v^{2},
v^{3}, v^{4})$ be two vectors in $\chi M$, then
$g(q^{3}u, q^{3}v)=g(q^{2}u, q^{2}v)=g(qu, qv)=g(u, v)$.
\end{prop}
\begin{thm}
Let $M$ be a Riemannian manifold with a metric $g$ from (\ref{f1})
and affinor structure from (\ref{f3}).  Let $\nabla$ be the
Riemannian connection of $g$. Then $\nabla q=0$ if and only if,
when
\begin{equation}\label{f4}
    \grad A=(\grad C)q^{2};\quad 2\grad B= (\grad C)(q+q^{3}).
\end{equation}
\end{thm}
\begin{proof}
Let $\Gamma_{ij}^{s}$ be the Christoffel symbols of $\nabla$. Let
$\nabla q=0$. That means
\begin{equation}\label{f5}
\nabla_{i}q^{s}_{j}=\partial_{i}q^{s}_{j}+\Gamma_{ik}^{s}q^{k}_{j}-\Gamma_{ij}^{k}q^{s}_{k}=0
\end{equation}
 From (\ref{f3}) and (\ref{f5}) we get
\begin{equation}\label{f6}
\Gamma_{ik}^{s}q^{k}_{j}=\Gamma_{ij}^{k}q^{s}_{k}
\end{equation}
Using (\ref{f1}), (\ref{f2}), (\ref{f3}), (\ref{2.1}), (\ref{f6})
and the well known identities:
\begin{equation}\label{2.3}
2\Gamma_{ij}^{s}=g^{as}(\partial_{i}g_{aj}+\partial_{j}g_{ai}-\partial_{a}g_{ij}).
\end{equation}
 After a long computation we get the following system:

 \begin{align*}
  &A_{4}-B_{1}+B_{3}-C_{2}=0,\\ &A_{4}+B_{1}-B_{3}-C_{2}=0,\\
 &2A_{2}+A_{4}-3B_{1}-B_{3}+C_{2}=0,\\
&A_{3}+B_{2}-B_{4}-C_{1}=0,\\
 &A_{3}-B_{2}+B_{4}-C_{1}=0,\\ &A_{2}-B_{1}+B_{3}-C_{4}=0,\\
 &A_{2}+B_{1}-B_{3}-C_{4}=0,\\ &A_{4}-B_{1}+3B_{3}+C_{2}+2C_{4}=0,\\ &A_{2}+2A_{4}-3B_{1}-B_{3}+C_{4}=0,\\
&A_{2}+2A_{4}-B_{1}-3B_{3}+C_{4}=0,\\ &A_{1}+2A_{3}-3B_{2}-B_{4}+C_{3}=0,\\ &A_{1}-B_{2}+B_{4}-C_{3}=0,\\
 &A_{3}-B_{2}-3B_{4}+C_{1}+2C_{3}=0,\\ &A_{1}-B_{2}-3B_{4}+2C_{1}+C_{3}=0,\\ &2A_{1}+A_{3}-B_{2}-3B_{4}+C_{1}=0,\\
 &A_{2}-B_{1}-3B_{3}+2C_{2}+C_{4}=0.
\end{align*}

The last system implies:
 \begin{align}\label{system2}\nonumber
   &A_{1}=C_{3},\  A_{2}=C_{4},\  A_{3}=C_{1},\  A_{4}=C_{2},\  B_{1}=B_{3}\\&B_{2}=B_{4},\  2B_{1}=C_{4}+C_{2},\ 2B_{2}=C_{1}+C_{3}.
\end{align}
 From (\ref{system2}) we find that (\ref{f4}) is valid.

Inversely, let (\ref{f4}) be valid. We can verify that (\ref{system2})
is valid, too. The identities (\ref{system2}) imply (\ref{f6}) and
consequently (\ref{f5}) is true. So $\nabla q=0$.
\end{proof}

Note. In fact (\ref{f4}) is a system of partial differential
equations for the functions $A$, $B$ and $C$. We can say that
(\ref{f4}) has a solution.

Now, we will give an example of such a manifold. Let
\begin{align}\label{example}\nonumber
    A&=(x^{1})^{2}+(x^{2})^{2}+(x^{3})^{2}+(x^{3})^{2},\\
    B&=x^{1}x^{2}+x^{2}x^{3}+x^{1}x^{4}+x^{3}x^{4},\\ \nonumber C&=2x^{1}x^{3}+2x^{2}x^{4}
\end{align}
 be three functions of a point $p(x^{1}, x^{2}, x^{3}, x^{4})\neq (x, x,
    x, x)$, $p(x^{1}, x^{2}, x^{3}, x^{4})\neq (-x, x,
    -x, x)$. Then $A > C > B > 0$ and
    \begin{equation}\label{fg}
    g_{ij}=\begin{pmatrix}
      A & B & C & B \\
      B & A & B & C \\
      C & B & A & B \\
      B & C & B & A \\
    \end{pmatrix}
\end{equation}
is positively defined. Also, we obtain $\grad A=(\grad
C)q^{2};\quad 2\grad B= (\grad C)(q+q^{3})$, which implies $\nabla
q=0$. So, we find an example for a manifold $M$ with a metric $g$, defined by
(\ref{example}) and (\ref{fg}), and affinor structures $q$,
defined by (\ref{f3}), which satisfies $\nabla q=0$.

Let $R$ be the curvature tensor field of $\nabla$, i.e $R(x,
y)z=\nabla_{x}\nabla_{y}z-\nabla_{[x,y]}z$. We consider the
associated with $R$ tensor field $R$ of type $(0, 4)$, defined by
the condition
\begin{equation*}
    R(x, y, z, u)=g(R(x, y)z,u), \qquad x, y, z, u\in \chi M.
    \end{equation*}
\begin{thm}
If $M$ is the Riemannian manifold with a metric $g$ and a parallel
structure $q$, defined by (\ref{f1}) and (\ref{f3}), respectively,
then the curvature tensor $R$ of $g$ satisfies the identity:
\begin{equation}\label{3.1}
    R(x, y, z, qu)=R(x, y, q^{3}z, u),\qquad x, y, z, u\in \chi M.
\end{equation}
\end{thm}
\begin{proof}
In the terms of the local coordinates $\nabla q=0$ implies
\begin{equation}\label{3.2}
    R^{l}_{sji}q_{k}^{.s}=R^{s}_{kji}q_{s}^{.l} .
\end{equation}
Using (\ref{f3}) we can verify
$q^{.i}_{j}=q_{a}^{.t}q_{j}^{.a}q_{t}^{.i}$ and then from
(\ref{f1}), (\ref{f3}) and (\ref{3.2}) we obtain (\ref{3.1}).
\end{proof}

\vspace{6pt}
\author{Dimitar Razpopov\\ Department of Mathematics and Physics\\ Agricultural University of Plovdiv\\
Bulgaria 4000\\
e-mail:dimitrerazpopov@hotmail.com}
\end{document}